\documentclass[11pt,twoside,a4paper]{article}
\usepackage{mathrsfs}
\usepackage{cite}
\usepackage{subfigure}
\usepackage{graphicx}
\DeclareFontFamily{U}{wncy}{}
\DeclareFontShape{U}{wncy}{m}{n}{<->wncyr10}{}
\DeclareSymbolFont{mcy}{U}{wncy}{m}{n}
\DeclareMathSymbol{\Sh}{\mathord}{mcy}{"58} 

\renewcommand{\vec}[1]{\mathbf{#1}}
\newtheorem{theorem}{Theorem}[section]

\newtheorem{lemma}[theorem]{Lemma}
\newtheorem{proposition}[theorem]{Proposition}
\newtheorem{corollary}[theorem]{Corollary}

\usepackage[cmex10]{amsmath}
\usepackage{amssymb}  
\usepackage{algorithmic}

\begin{document}
\title{
On the Discrepancy Normed Space of Event Sequences for Threshold-based Sampling
}

\author{Bernhard A. Moser\\
Software Competence Center Hagenberg, Austria\\
Email: bernhard.moser@scch.at }

\maketitle

\begin{abstract}
Recalling recent results on the characterization of threshold-based sampling as quasi-isometric mapping, mathematical implications on the metric and topological structure of the space of event sequences  are derived. In this context, the space of event sequences is extended to a normed space equipped with Hermann Weyl's discrepancy measure. Sequences of finite discrepancy norm are characterized by a Jordan decomposition property. Its dual norm turns out to be the norm of total variation. As a by-product a measure for the lack of monotonicity of sequences is obtained. A further result refers to an inequality between the discrepancy norm and total variation which resembles Heisenberg’s uncertainty relation.  
\end{abstract}


\section{Motivation}
\label{s:motivation}
This paper starts by recalling a recent result for the understanding of threshold-based sampling schemes as quasi-isometric mapping~\cite{Moser17}.
In this context a threshold-based sampling scheme is understood as a mapping from the space of sampled signals to the space of resulting {\it event sequences} of ``up'' and ``down'' events that preserves the notion of ``closeness'' or synonymously ``similarity''. 
The ``up'' and ``down'' events are triggered by the sampling process. Usually these events are represented by $+1$ and $-1$, respectively.
To be precise, preserving the topology is not possible in the strict sense (see e.g.~\cite{MoserEBCCSP16}). 
This effect is an immediate consequence of the all-or-nothing law of threshold-based sampling.
Either there is a triggering sampling event at a certain time or not. 
Take for example signals below threshold. Such signals cannot be distinguished from the samples, because there are none. 
So preserving the metric as e.g. the notion of closeness can only be satisfied in a relaxed fashion, namely as quasi-isometry.
As a consequence, we single out metrics being compatible with the quasi-isometry constraint. 
As pointed out in~\cite{Moser17}, this analysis leads to the class of metrics for which a sequence of alternating ``up'' (modeled by $1$) and ``down'' (modeled by $-1$) events is considered to be close to the zero sequence that contains no event at all. 

One metric that fulfils this condition is due to Hermann Weyl, namely the so-called discrepancy measure~(see,~\cite{Weyl1916,Doerr2014, Moser2012a}). This measure was introduced over 100 years ago in the context of evaluating the quality of pseudo-random numbers. In a vector space this measure leads to a norm, the discrepancy norm $\|.\|_{\mbox{\tiny D}}$. This norm distinguishes itself from the familiar Euclidean or another $L_p$ norm by its asymmetric shape of its unit ball. This asymmetry is due to the fact that the norm evokes in general different lengths after rearranging the order of events in a sequence.
There is an instructive interpretation of the discrepancy. Consider a walker along a line, who makes a step ahead if the event is ``up''
and a step backwards, if the event is ``down''. The discrepancy is the range of the walk.

As shown in~\cite{MoserTSP2014}, typical metrics in this context such as the van Rossum~\cite{Rossum01}  or the Victor-Purpura metric~\cite{Victor1996} do not satisfy this condition. As a consequence arbitrary small deviations can cause disruptive effects in the input-output behavior when relying on similarity measures based on such metrics.


In this paper we focus on mathematical implications on the  topological structure of the space of event sequences when underlying Weyl's discrepancy norm $\|.\|_{\mbox{\tiny D}}$. As first result, we provide a characterization of those event sequences that are finite in this metric in a way that resembles the Jordan decomposition law of functions of total variation, see Section~\ref{s:Jordan}.  
This result indicates a close relationship between the discrepancy norm, $\|.\|_{\mbox{\tiny D}}$, 
and the semi-norm of total variation, $\|.\|_{\mbox{\tiny BV}}$.
In analogy to $L_p$ spaces we denote the space of event sequences that are bounded with respect to $\|.\|_{\mbox{\tiny D}}$ 
 by $L_{\mbox{\tiny D}}$.
In Section~\ref{s:dual} we study the dual space $L^*_{\mbox{\tiny D}}$ of $L_{\mbox{\tiny D}}$.
As second result we identify $L^*_{\mbox{\tiny D}}$ as the space of functions of total variation.

As measure of oscillation  $\|.\|_{\mbox{\tiny BV}}$ behaves inverse proportional to $\|.\|_{\mbox{\tiny D}}$.
If the range of a walk consisting of $+1$ and $-1$ steps of length $n \in \mathbb{N}$ is small then there is much oscillation. 
For example, for a sequence of alternating signs, $+1, -1, +1 \ldots$, the range is minimal and the oscillation is maximal.
On the other hand, little oscillation means that there is a predominant direction of the walk and therefore a larger range.
This reciprocal relation is topic of Section~\ref{s:inequality}, which leads to the inequality ($x_i \in \{-1,1\}$, 
$\vec x = (x_1, \ldots, x_n)$ not constant)
$$
n \leq \|\vec x\|_{\mbox{\tiny D}} \cdot \|\vec x\|_{\mbox{\tiny BV}},
$$ 
which resembles Heisenberg's uncertainty relation in its form. 
On the left hand side there is a constant as lower bound and on the right hand side there is a product of two measures that represent dual concepts.
For the Heisenberg inequality these dual concepts are time and frequency. 
In our case, the dual concepts refer to oscillations in terms of total variation and range of the corresponding walk.

Before, we start with a section on preliminaries (Section~\ref{s:intro}) by introducing and recalling the notion of quasi-isometry 
(Subsection~\ref{ss:distances}), Weyl's discrepancy (Subsection~\ref{ss:discrepancy}) and its relation to quasi-isometry in the context of threshold-based sampling (Subsection~\ref{ss:quasi}).

\section{Preliminaries}
\label{s:intro}
\subsection{Mathematics of Distances}
\label{ss:distances}
First of all, let us fix some notation. $1_{I}$ denotes the indicator function of the set $I$, i.e., $1_I(t)=1$ if $t\in I$ and $1_I(t)=0$ otherwise.  
 $\|.\|_{\infty}$ denotes the uniform norm, i.e., $\|f-g\|_{\infty} = \sup_{t \in X} |f(t)-g(t)|$, where $X$ is the domain 
of $f$ and $g$. If $M$ is a discrete set then $|M|$ denotes its number of elements. If $I$ is  an interval, then 
$|I|$ denotes its length. $\mathcal{I}$ denotes the family of real intervals.

In this section we recall basic notions related to distances such as semi-metric, isometry and quasi-isometry, see
 e.g.,~\cite{EncyclopediaofDistances2009}.

Let $X$ be a set. A pseudo-metric $d:X\times X \rightarrow [0,\infty)$ is characterized by a) $d(x,x)=0$  for all $x\in X$, b) $d(x,y)=d(y,x)$ for all $x,y \in X$ and c) the triangle inequality
$d(x,z)\leq d(x,y) + d(y,z)$ for all $x,y,z \in X$. 
$d$ is a metric if, in addition to a) the stronger condition a')
$d(x,y)=0$ if and only if $x=y$, is satisfied.
The semi-metric  $\tilde d$ is called {\it equivalent} to $d$, in symbols $d \sim \tilde d$, if and only if there are constants 
$A_1, A_2>0$ such that
\begin{equation}
\label{eq:norm-equivalencecondition}
A_1 d(x,y)  \leq \tilde d(x,y) \leq A_2 \, d(x,y)
\end{equation}
 for all $x$, $y$ of the universe of discourse.

A map $\Phi: X\rightarrow Y$ between a metric space ${(X,d_{X})}$ and another metric space  
$(Y,d_{Y})$ is called {\it isometry} if this mapping is distance preserving, i.e., for any 
$x_1,x_2 \in X$ we have $d_{X}(x_1,x_2) = d_{Y}(\Phi(x_1), \Phi(x_2))$.

The concept of {\it quasi-isometry}  relaxes the notion of  isometry by imposing only a 
coarse Lipschitz continuity and a coarse surjective property of the mapping. 
$\Phi$ is called a {\it quasi-isometry} from 
$(X,d_{X})$ to $(Y,d_{Y})$ if there exist constants 
$A\geq 1$, $B\geq 0$, and $C\geq 0$ such that the following two properties hold:
\\
\noindent
i) For every two elements $x_1, x_2\in X$, the distance between their images is, up to the additive constant $B$, within a factor of $A$ of their original distance. This means,  there are constants $A$ and $B$ such that 
$\forall x_1, x_2\in X$
\begin{equation}
\label{eq:quasi1}
{\frac{1}{A}}\,d_{X}(x_1,x_2)-B 
 \leq d_{Y}(\Phi(x_1),\Phi(x_2))  
 \leq A\,d_{X}(x_1,x_2)+B. 
\end{equation}
\\
\noindent
ii) Every element of $Y$ is within the constant distance $C$ of an image point, i.e.,
\begin{equation}
\label{eq:quasi2}
\forall y \in Y:\exists x\in X:d_{Y}(y,\Phi(x))\leq C.
\end{equation}

Note that for $B=0$ the condition (\ref{eq:quasi1}) reads as Lipschitz continuity condition of the operator $\Phi$.
This means that (\ref{eq:quasi1}) can be interpreted as a relaxed bi-Lipschitz condition.
The two metric spaces $ (X,d_{X})$ and $(Y,d_{Y})$ are called {\it quasi-isometric} if there exists a quasi-isometry $Q$ from $ (X,d_{X})$ to $ (Y,d_{Y})$.
  
In this paper, the total variation $\|.\|_{\tiny \mbox{BV}}$ plays a central role.
It is a measure for the amount of oscillations and is defined by
\begin{equation}
\label{eq:BV}
\|f\|_{\tiny \mbox{BV}}:= \sup_{x_1 < x_2 \ldots < x_n | x_i \in \mathbb{R}, n \in \mathbb{N}} \sum_{i=1}^{n-1} | f(x_{i+1}) - f(x_i)|.
\end{equation}
(\ref{eq:BV}) is a semi-norm as $\|f\|_V=0$ does not imply $f=0$, rather any constant function has a vanishing total variation.
For a sequence $\vec x = (x_i)_i$ the total variation equals $\sum_i |x_{i+1}-x_i|$.

\subsection{Introduction to Weyl's Discrepancy}	
\label{ss:discrepancy}	
\begin{figure}[htp]
  \begin{center}   	      
	    \includegraphics[width=0.80 \columnwidth]{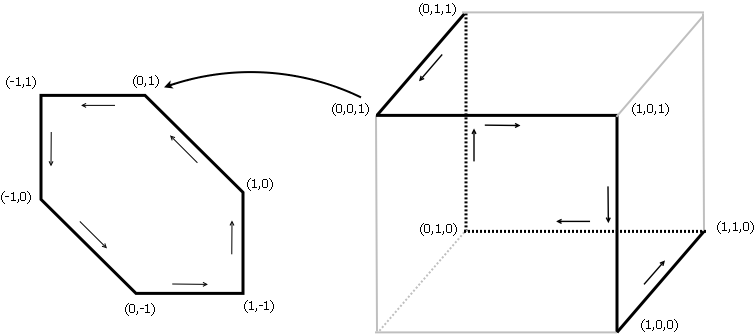} 
  \end{center}
  \caption{The vertices from the hypercube (right) are mapped by SOD with threshold $\theta =1$ to vertices of the unit ball of the discrepancy norm (left)}
  \label{fig:SODHypercube}
\end{figure}

Let us make a thought experiment by considering a discrete sequences of $0$ and $1$ entries, i.e.,
$\vec x = (x_1, \ldots, x_n) \in \{0,1\}^n$. This sequence actually is a vertex in the $n$-dimensional hypercube
$[0,1]^n$.
Now, let us apply the Send-on-Delta (SOD) sampling scheme with threshold $\theta=1$. For example, $(0,1,0)$ is mapped to $(1,1)$, as $x_2-x_1\geq 1$ yields $1$ and 
$x_3-x_2\leq -1$ yields $-1$. How does this operation transform the geometry of the hypercube?
See Figure~\ref{fig:SODHypercube} for an illustration of the resulting polytope for $n=3$.
The resulting body is a bounded convex polytope $P$ that is point-symmetrical w.r.t the origin. 
Such a body $P$ constitutes a norm 
$$
\|\vec x \|:= \left( 
\max\{ \lambda>0 \mid \lambda \vec x \in P\}
\right)^{-1}.
$$
As shown in~\cite{Moser2012a} the norm induced by SOD on the hypercube yields Hermann Weyl's discrepancy.

In~\cite{Weyl1916} Weyl  introduces a concept of discrepancy in the context of
  pseudo-randomness of sequences of numbers from the unit interval.
Weyl's discrepancy concept leads to the definition
\begin{equation}
\label{discrepancy_vector}
\|\vec x\|_D =  \sup_{n_1, n_2 \in \mathbb{Z}: n_1 \leq n_2,} |\sum_{i=n_1}^{n_2} x_i|,
\end{equation}
which induces a norm on the $n$-dimensional real vector space~\cite{Moser2012a}.
Applications of the norm~(\ref{discrepancy_vector}) can be found in pattern recognition~\cite{neunzertWetton87}, print inspection in the context of pixel classification~\cite{BauerBodenhoferKlement96}, template matching and registration~\cite{Moser2008}.
In contrast to $p$-norms $\|.\|_p$, $\|\vec x\|_p = (\sum_i |x_i|^p)^{(1/p)}$,
the norm $\|.\|_D$ strongly depends on the sign and also the ordering of the entries, as illustrated by the examples  
$\|(-1,1,-1,1)\|_D = 1$ and  $\|(-1,-1,1,1)\|_D = 2$.

\begin{figure}[htp]
  \begin{center}
   \subfigure[]
	      {\label{fig:monotonieLipschitz:n2fmax}
	      \includegraphics[width=0.45 \columnwidth]{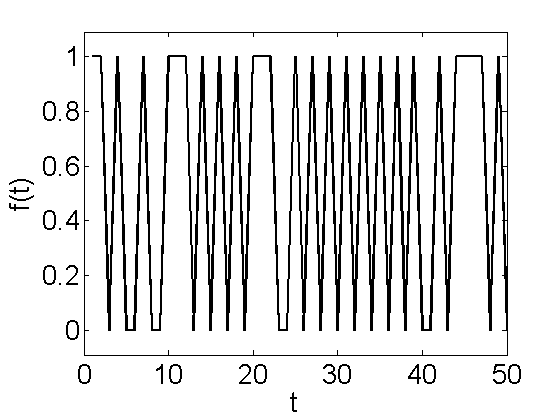}} 
    \subfigure[]
	      {\label{fig:monotonieLipschitz:dn}
	      \includegraphics[width=0.45 \columnwidth]{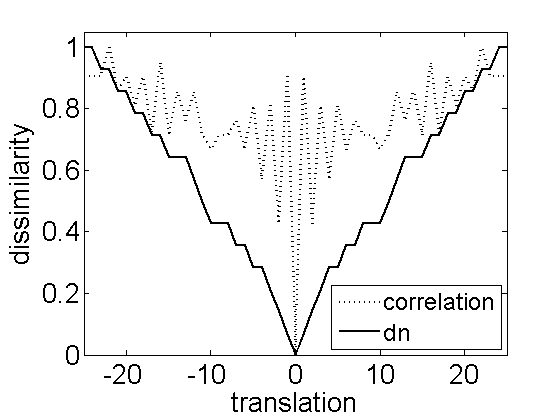}} 
  \end{center}
  \caption{Figure (a) shows a sequence of all-or-none events. Figure (b) depicts its misalignment function with respect to the Euclidean norm (dashed line) and with respect to the discrepancy norm (solid line). Note that the solid line is monotonic according to (P6)}.
  \label{fig:monotonieLipschitz}
\end{figure}

Generally, $\vec x = (x_i)_i$ with $x_i\geq 0$ entails $\|\vec x\|_D = \|\vec x\|_{1}$, and 
$\vec x = ((-1)^{i})_i$ the equality $\|\vec x\|_D = \|\vec x\|_{\infty}$, respectively, indicating that the more there are alternating signs of consecutive entries, the lower is the value of the discrepancy norm. 
Observe that $
 \|\vec x\|_{\infty} \leq \|\vec x\|_D \leq \|\vec x\|_1, 
$
hence, due to Hoelder's inequality
$
 n^{-1/p} \|\vec x\|_p \leq \|\vec x\|_D \leq  n^{1-1/p} \|\vec x\|_p
$.
For convenience let us consider a sequence $(x_i)_i$ with $i \in I_n$,
$x_i=0$ for $i \notin I_n$, and denote by $\Delta_{\vec x}(k)= \|(x_{i+k} - x_{i})_i\|_D$ the misalignement function of $x$ with respect to $\|.\|_D$.
Then we have the following properties~\cite{Moser2008}:
\begin{description}
\item[(P1)] $\|(x_i)_{i\in I_n}\|_D$ induces a norm on $\mathbb{R}^n$. 
\item[(P2)] $\Delta_{\vec x}(0)=0$ for all summable real sequences $\vec x$.
\item[(P3)] 
$\|
(x_i)_{i\in I_n}\|_D = 
\max\{0,\max_{k\in I_n} \sum_{i = 1}^k x_i\} 
- \min\{0,	\min_{k \in I_n} \sum_{i = 1}^k x_i\}$ 
\item[(P4)] Lipschitz property: $\Delta_{\vec x}(k) \leq |k| \cdot L$,   where  
$L = \max_i x_i - \min_i x_i$ and $k \in \mathbb{Z}$. 
\item[(P5)] $\Delta_{\vec x}(k) = \Delta_{\vec x}(-k)$ for 
$\vec x= (x_i)_i$ with $x_i \geq 0$ and $k \in \mathbb{Z}$. 
\item[(P6)] For $\vec x= (x_i)_i$ with $x_i \geq 0$ the function $\Delta_{\vec x}(.)$ is monotonically increasing on $\mathbb{N}\cup \{0\}$.
\end{description}	
Equation (P3) allows us to compute the discrepancy of a sequence of length $n$ with $O(n)$ operations instead of $O(n^2)$ number of operations resulting from the original Definition~(\ref{discrepancy_vector}).
Especially the monotonicity (P6) as well as the Lipschitz property (P4) are interesting properties for
applications in the field of signal analysis. It is instructive to point out that the Lipschitz constant in (P4) does not depend on frequencies or other characteristics of the sequence $\vec x$. 
Properties (P4), (P5) and (P6) are illustrated in the Figures~\ref{fig:monotonieLipschitz:n2fmax} 
and~\ref{fig:monotonieLipschitz:dn} which demonstrate the behavior of the misalignment function of a sequence of all-or-none events. While Figure~\ref{fig:monotonieLipschitz:n2fmax} shows typical local minima of the misalignment function
 with respect to the Euclidean norm, Figure~\ref{fig:monotonieLipschitz:dn} visualizes the symmetry property (P5),  the monotonicity property (P6) and the boundedness of its slope due to the Lipschitz property (P4) of the corresponding misalignment function induced by the discrepancy norm. 

To obtain a clear interpretation of the discrepancy, let's think of a walker who moves up or down along a line at each time step  
according to the sequence $(x_1, \ldots,x_n) \in \{-1,1\}^n$. What is the range of this movement?
Consider the pair of variables $(t,d)$ for time and distance. The walk can be represented by the graph
$$
\gamma = \left((0,0)^T, (1,x_1)^T, \ldots, \sum_{i=1}^n (1,x_i)^T \right)
$$
in $(\mathbb{N}_0 \times \mathbb{Z})^{n+1}$.
The diameter (range) of $\gamma$ w.r.t. the second variable, i.e., in the direction of $(0,1)^T$,  is given by
\begin{eqnarray}
\label{eq:disc}
& & \max_{1 \leq n_1, n_2 \leq n} 
\left| \sum_{i=n_1}^{n_2} \left\langle (1,x_i)^T, (0,1)^T\right\rangle \right| \nonumber \\
& = &  
\max_{1 \leq i\leq n} \sum_{j=0}^i x_j - \min_{1 \leq i\leq n} \sum_{j=0}^i x_j \nonumber \\
& = &
\|\vec x\|_{\tiny \mbox{D}}
\end{eqnarray}
where $\left\langle.,.\right\rangle$ denotes the usual inner product, and $x_0 = 0$.
Equation~(\ref{eq:disc}) tells us that the discrepancy can be interpreted as range. 
It is interesting to note that this interpretation was the key to solve the 
problem of computing the distribution of the range of a random walk~\cite{Moser2014Random}.
A problem that remained unsolved for more than 50 years after it was stated by Feller in 1951~\cite{fel51}.

\subsection{Threshold-Based Sampling as Quasi-Isometry }	
\label{ss:quasi}	
\cite{Moser17} provides a framework for constructing metrics in the input and the output space of a threshold-based sampling scheme
$\Phi_{\theta}$ such that $\Phi_{\theta}$ becomes a quasi-isometry with constants 
$A_{\theta}$ and $B_{\theta}$, according to (\ref{eq:quasi1}). 
The construction relies on Weyl's discrepancy norm. The metrics can be constructed in a way that $\lim_{\theta \rightarrow 0} A_{\theta} = 1$ and $\lim_{\theta \rightarrow 0} B_{\theta} = 0$ (see Theorem 6.1 of \cite{Moser17}). This means that these metrics are asymptotically isometric for ever decreasing thresholds.

For example, for Send-on-Delta (SOD) and Integrate-and-Fire (IF) we obtain $A_{\theta} = 1$ and $B_{\theta} = 4 \theta$. 
In both cases we obtain $\|.\|_{\tiny \mbox{D}}$ as metric in the output space, that is the space of event sequences.
In the input space, in the former case (SOD) we get the semi-norm of the range and for the latter (IF) we obtain as metric an integral version 
of the discrepancy norm for integrable functions.

Further analysis shows that the choice of the discrepancy measure or some quasi-isometric variant of it is even necessary in order to
turn $\Phi_{\theta}$ into a quasi-isometry. 
	
This special role of the discrepancy measure in the context of threshold-based sampling 
strongly motivates to investigate the space of event sequences based on the discrepancy measure as metric in more detail.

\section{Conception of the Space of Event Sequences as Metric Space}
\label{s:Jordan}
Taking up the results about quasi-isometry of Subsection~\ref{ss:quasi}, we come up with the following postulates 
for the space of event sequences for threshold-based sampling.

Basically, an event sequence is a function in time that is zero except at discrete time points of triggered events.
In the case of bipolar events we therefore have functions of the form $\eta:[0,\infty) \rightarrow \{-1,0,1\}$.

As the events are triggered by the sampling scheme, the events are sparse, that is there are no accumulation points of events.
Putting in other words, for any finite time interval $[a,b]$ there are only a finite number of events inside this interval.

Now, let us extend this space to the vector space of functions $\eta:[0,\infty) \rightarrow \mathbb{Z}$ and 
equip this space with the discrepancy norm $\|.\|_{\tiny \mbox{D}}$. 
Note that an event sequence can synonymously be represented by its sequence of events $(t_k, \eta_k)_k$ ($\eta_k:=\eta(t_k)$) 
which justifies the term ``sequence'' in this context.
Therefore, the discrepancy norm of an event sequence, $\|\eta\|_{\tiny \mbox{D}}$ is well defined by referring to the sequence, i.e.,
\begin{eqnarray}
\label{eq:norm_e}
\|\eta\|_{\tiny \mbox{D}} & := &\|(\eta_k)_k\|_{\tiny \mbox{D}} \nonumber\\
& =  & \sup_{[a,b]} \left|\int_a^b \eta dc\right| = \sup_{[a,b]} \left|\sum_{k = a}^b \eta(t_k) \right|,
\end{eqnarray}
where the last line of (\ref{eq:norm_e}) represents the sum as integral w.r.t. the counting measure $c$.

Let us denote this normed space of event sequences $\eta:[0,\infty) \rightarrow \mathbb{Z}$ of finite discrepancy, 
$\|\eta\|_{\tiny \mbox{D}}< \infty$, by 
\begin{equation}
\label{eq:space_e}
(\mathscr{E}_{\tiny \mbox{D}}, \|.\|_{\tiny \mbox{D}}).
\end{equation}

Analogously, referring to the input space of signals we can equip the the space of locally integrable functions, $\mathscr{P}(\mathbb{R})$, with the discrepancy in its integral version 
\[
\|f\|_{\tiny \mbox{D}, \lambda}  := \sup_{[a,b]} \left|\int_a^b f d\lambda\right|
\]
w.r.t. the Lebesgue measure $\lambda$.
Let us denote
\begin{eqnarray}
\label{eq:norm_f}
\mathscr{L}_{\tiny \mbox{D}} &:=&  \{f \in \mathscr{P}(\mathbb{R})\mid \|f\|_{\tiny \mbox{D}, \lambda} < \infty\}.
\end{eqnarray}
We refer to the corresponding normed space by 
\begin{equation}
\label{eq:space_f}
(\mathscr{L}_{\tiny \mbox{D}}, \|.\|_{\tiny \mbox{D}, \lambda}).
\end{equation}

Next, we present the results of this paper.
First we provide a characterization  of event sequences and functions of finite discrepancy.

\subsection{Jordan Decomposition of Finite Discrepancy Sequences and Functions}	
\label{s:Jordan}
The following Lemma~\ref{lem:diam} shows that $\|f\|_{\tiny \mbox{D}, \lambda}$ can be represented as 
range of values assumed by the function 
\begin{eqnarray}
\label{eq:Gamma}
\Gamma_f(x) &:=& \liminf_{n \in \mathbb{Z}, n < x} \int_{-n}^x f d\lambda, \\
     & = & \liminf_{n \rightarrow \infty} \int_{-n}^x f d\lambda. \nonumber 
\end{eqnarray}
\begin{lemma}({\bf Discrepancy as Range, $\lambda$-Version})
\label{lem:diam}

Let $f \in \mathscr{P}(\mathbb{R})$. Then,
\begin{equation}
\label{eq:diam}
\|f\|_{\tiny \mbox{D}, \lambda}   = 
\sup_x  \Gamma_f(x) - \inf_x \Gamma_f(x). 
\end{equation}
\end{lemma}

First of all note that the assumption 
$$\|f\|_{\tiny \mbox{D}, \lambda} = \sup_{a,b}\left|\int_a^b f d\lambda \right|< \infty$$
guarantees that $\liminf_{n \rightarrow \infty} \int_{-n}^x f d\lambda$ is well defined and finite for all 
real $x$.

Now, observe that the assumption
$$ \inf_x \Gamma_f(x) = \varepsilon >0$$
 implies
that for all $k \in \mathbb{Z}$ there is a natural number $N_k\in \mathbb{N}$ such that for 
all $n \geq N_k$ there holds $\inf_{m\geq n} \int_{-m}^k f d\lambda \geq \varepsilon/2$. Hence, there is a sequence of increasing numbers $(k_n)_n$, $\lim_n k_n = \infty$, such that  
$\int_{-k_{n+1}}^{k_n} f d\lambda \geq \varepsilon/2$.
Consequently, we obtain $\sup_{a,b}|\int_a^b f d\lambda | \geq \int_{-k_{n+1}}^{k_1} f d\lambda \geq n\,\varepsilon/2$ which contradicts 
$\|f\|_{\tiny \mbox{D}, \lambda}  < \infty$. Hence,  
\begin{equation}
\label{eq:c}
\inf_x \Gamma_f(x)   \leq 0.
\end{equation}
Analogously, we obtain  
\begin{equation}
\label{eq:d}
\sup_x \Gamma_f(x) \geq 0.
\end{equation}

Further, note that~(\ref{eq:diam}) yields
\begin{equation}
\label{eq:split}
\int_a^b f d\mu = \liminf_{n \rightarrow \infty} \int_{-n}^b f d\lambda- \liminf_{n \rightarrow \infty} \int_{-n}^a f d\lambda
\end{equation}
for $a < b$.
Taking~(\ref{eq:c}),~(\ref{eq:d}) and~(\ref{eq:split}) together proves Lemma~\ref{lem:diam}.
$\Box$
For an example take $f(t) = \sin(t)$ on $\mathbb{R}$. 
This function has finite discrepancy, namely $\|f\|_{\tiny \mbox{D}} = \int_0^{\pi} \sin(t)dt= 2$. 
Note that in general locally integrable periodic functions have finite discrepancy. 

In analogy to (\ref{eq:Gamma}) we define
\begin{eqnarray}
\label{eq:gamma}
\gamma_{\eta}(k) &:=& \liminf_{n \rightarrow \infty} \sum_{j= -n}^k \eta_j. \nonumber 
\end{eqnarray}
and obtain an identity in analogy to Lemma~\ref{lem:diam}, i.e.,
\begin{equation}
\label{eq:diamDiscrete}
\|(\eta_k)_k\|_{\tiny \mbox{D}} = \sup_{k\in \mathbb{Z}} \gamma_{\eta}(k) - \inf_{k\in \mathbb{Z}} \gamma_{\eta}(k), 
\end{equation}
where $(\eta_k)_k \in {\mathbb{R}}^{\mathbb{Z}}$.

It is a well known result, the so-called Jordan decomposition law, that functions $f$ of bounded variation can be characterized  as difference
of monotonic functions $h_1$ and $h_2$, $f = h_2 - h_1$. 
For functions of bounded discrepancy we obtain an analogous characterization.

\begin{theorem}({\bf Jordan Decomposition of Bounded Discrepancy, $\lambda$-Version})
\label{cor:mon}

Let $f \in \mathscr{P}([a,b])$, $a< b$. Then 
$ \|f\|_{\tiny \mbox{D}, \lambda} \leq r < \infty$ if and only if
there are non-decreasing locally absolutely continuous functions $h_1$, $h_2$  such that
$\|h_2 - h_1\|_{\infty} \leq r/2$ and  $f = \dot{h_2} - \dot{h_1}$ almost everywhere.
\end{theorem}

For the proof we split $f = f^+ - f^-$ into its non-negative and non-positive part
$f^+ = \max\{f, 0\}$,  $f^- = \min\{f, 0\}$.
For $\|f\|_{\tiny \mbox{D}, \lambda} = 0$ we choose $h_1= h_2 = 0$. 

Further on, assume that $\|f\|_{\tiny \mbox{D}, \lambda} >  0$.
Due to the compactness of $[a,b]$ there is an interval $[a^*, b^*] \subseteq [a,b]$ such that
$r:=\|f\|_{\tiny \mbox{D}, \lambda} = |\int_{a^*}^{b^*} f d\lambda|$. 
Due to the intermediate value theorem there is a real $c^* \in [a^*, b^*] $ such that
$|\int_{a^*}^{\overline{x}} f d\lambda| = |\int_{c^*}^{b^*} f d\lambda| = r/2$. 
Let us define
\begin{eqnarray}
h_2(x) := \int_{c^*}^x f^+ d \lambda, \,\, h_1(x) := -\int_{c^*}^x f^- d \lambda. \nonumber
\end{eqnarray}
Consider the intervals $[a_k, b_k]$, $k \in \mathbb{Z}$,   at which $|\int_{a_k}^{b_k} f d\lambda|$ assumes its maximum, that is
$|\int_{a_k}^{b_k} f d\lambda| = \|f\|_{\tiny \mbox{D}, \lambda} = r$. 
Note that $h_1$ and $h_2$ are non-decreasing and almost everywhere differentiable with $f = \dot{h_2} - \dot{h_1}$. 
Further, note that $\int_{a_k}^{b_k} f d\lambda \in \{-r, r\}$ is an alternating sequence which implies 
 that $|h_2(x) - h_1(x)| = |\int_{c^*}^{x} f d\lambda| \leq r/2$ for all $x$, hence $\|h_2 - h_1 \|_{\infty} \leq r/2$.

On the other hand, let us suppose that 
$\|h_2 - h_1 \|_{\infty} \leq r/2$ where $h_1, h_2$ are absolutely continuous functions satisfying
$f = \dot{h_2} - \dot{h_1}$ with $\dot{h_2}, \dot{h_1} \geq 0$ a.e.. Then  
Lemma~\ref{lem:diam}  entails
\begin{eqnarray}
\label{eq:proofJordan}
 & & \|f\|_{\tiny \mbox{D}, \lambda}   \\
																& = & \sup_x \Gamma_f(x) - \inf_x \Gamma_f(x)  \nonumber \\
																 & = & \sup_x \liminf_{n\rightarrow \infty} \int_{-n}^x \dot{h_2}(x) - \dot{h_1}(x) d\lambda  \nonumber\\
																 &   &  - \inf_x \liminf_{n\rightarrow \infty} \int_{-n}^x \dot{h_2}(x) - \dot{h_1}(x) d\lambda \nonumber\\
																 & = &  \sup_x (h_2(x) - h_1(x)) - \liminf_{n\rightarrow \infty}  (h_2(-n) - h_1(-n)) \nonumber\\
																 &   & - \inf_x (h_2(x) - h_1(x) ) + \liminf_{n\rightarrow \infty}  (h_2(-n) - h_1(-n)). \nonumber
\end{eqnarray}
Since $\|h_2 - h_1\|_{\infty} < \infty$ implies that  $\liminf_n (h_2(-n) - h_1(-n))$ exists, that is 
$\liminf_n (h_2(-n) - h_1(-n)) = \rho \in \mathbb{R}$, Equation~(\ref{eq:proofJordan}) finally implies
$$
\|f\|_{\tiny \mbox{D}, \lambda} =  \sup_x (h_2(x) - h_1(x)) - \inf_x (h_2(x) - h_1(x) ) \leq 2 \frac{r}{2}<\infty,
$$
which ends the proof.
$\Box$

In an analogous way we obtain a Jordan decomposition representation for the discrete version.
\begin{theorem}({\bf Jordan Decomposition of Bounded Discrepancy, Discrete Version})
\label{cor:mon}
Let $\eta = (\eta_k)_{k \in \mathbb{N}} \in {\mathbb{R}}^{\mathbb{N}}$. 
Then 
$ \|\eta\|_{\tiny \mbox{D}} = r < \infty$ if and only if
there are non-decreasing sequences 
$$
\chi_1 = (\chi_1(k))_k,\,\, \chi_2 = (\chi_2(k))_k
$$  
such that
$\|\chi_2 - \chi_1\|_{\infty} \leq r/2$ and  
$$
\eta(k) = (\chi_2(k) - \chi_2(k-1)) - (\chi_1(k) - \chi_1(k-1)).
$$
\end{theorem}

Assume that $ \|\eta\|_{\tiny \mbox{D}} = r < \infty$.
We set 
\begin{eqnarray}
\label{eq:disJordanChi}
\chi_2^{(\alpha)}(k) &:=& \sum_{i = 1}^k \max\{0, \eta_i\} - \alpha,  \\
\chi_1(k) &:= & -\sum_{i = 1}^k \min\{0, \eta_i\}, \nonumber
\end{eqnarray}
where $\alpha$ in the first line in (\ref{eq:disJordanChi}) is defined by
\[
\alpha:= \frac{1}{2}\left( 
\max_{k \in \mathbb{Z}}\{\chi_2^{(0)}(k) - \chi_1(k) \} - \min_{k \in \mathbb{Z}}\{\chi_2^{(0)}(k) - \chi_1(k)\}
\right).
\]
For convenience we define $\chi_2(0):= - \alpha$ and $\chi_1(0):=0$.

Note that $\chi_1$ and $\chi_2$ are non-decreasing. Further, we check that
$$
(\chi_2(k) - \chi_2(k-1)) - (\chi_1(k) - \chi_1(k-1))= \eta_k
$$ 
and that
$$
|\chi_2(k) - \chi_1(k) |  \leq r/2.
$$

The other direction of the proof follows.
Suppose $\eta(k) = (\chi_2(k) - \chi_2(k-1)) - (\chi_1(k) - \chi_1(k-1))$,
$\|\chi_2 - \chi_1\|_{\infty} \leq r/2$ 
and consider the range representation of the discrepancy
\begin{eqnarray}
\|\eta\|_{\tiny \mbox{D}} & = & \max_{k \in \mathbb{N}}\{0, \sum_{j=1}^k (\chi_2(j) - \chi_2(j-1))  \nonumber \\
													&    & - (\chi_1(j) - \chi_1(j-1)) \} \nonumber \\
													&    & - \min_{k \in \mathbb{N}}\{0, \sum_{j=1}^k (\chi_2(j) - \chi_2(j-1))  \nonumber \\
													&    & - (\chi_1(j) - \chi_1(j-1)) \} \nonumber \\
													& = & \max_{k \in \mathbb{N}}\{0, \chi_2(k) - \chi_1(k) \} \nonumber \\
													&    & - \min_{k \in \mathbb{N}}\{0, \chi_2(k) - \chi_1(k) \} \nonumber \\
													& \leq & r. \nonumber													
\end{eqnarray}
$\Box$

As a corollary we obtain the result that a bounded discrepancy function can also be characterized by a differentiable function whose range is bounded.
It turns out that the Lebesgue measure of this range equals the discrepancy.

\begin{corollary} ({\bf Discrepancy as Range, Second Version})
\label{th:characterization}
Let $f \in \mathscr{P}(\mathbb{R})$. Then 
$\|f\|_{\tiny \mbox{D}, \lambda} = r < \infty$ if and only if
there is a uniquely determined locally absolutely continuous  function 
$g$ such that $\overline{g(\mathbb{R})}=[0,r]$ and $f = \dot{g}$ almost everywhere.
\end{corollary}

Suppose $\|f\|_{\tiny \mbox{D}, \lambda} = r < \infty $.
Let us introduce 
\begin{eqnarray}
\label{eq:g}
g(x) & = & -c + \liminf_{n \rightarrow \infty}\int_{-n}^x f d\lambda 
\end{eqnarray}
where 
$ c  :=  \inf_{x} \liminf_{n \rightarrow \infty} \int_{-n}^x f d\lambda$.
Due to the fundamental theorem of Lebesgue integral calculus 
$g$ is locally absolutely continuous and differentiable almost everywhere. 
Equation~(\ref{eq:g}) implies $\inf_x g(x) = 0$ by construction.
Now, consider
\[
\sup_x g(x) =  - \inf_{x} \liminf_{n \rightarrow \infty} \int_{-n}^x f d\lambda +
\sup_x \liminf_{n \rightarrow \infty} \int_{-n}^x f d\lambda
\]
which by Lemma~\ref{lem:diam} yields
$\sup_x g(x) = \|f\|_{\tiny \mbox{D}, \lambda} = r$.
The identity $f = \dot{g}$ almost everywhere follows from  construction (\ref{eq:g}).

Now, consider an absolutely continuous function $g$ with 
$\overline{g(\mathbb{R})}=[0,r]$, i.e., 
$\inf_x g(x) =0$ and $\sup_x g(x) = r \geq 0$.
 
Then,  
\begin{eqnarray}
\sup_x  \liminf_{n \rightarrow \infty} \int_{-n}^x \dot{g} d\lambda & = & 
\sup_x  \liminf_{n \rightarrow \infty} (g(x) - g(-n)) \nonumber\\
& = & 
\sup_x  g(x) - \limsup_{n \rightarrow \infty} g(-n) \nonumber
\end{eqnarray}
and, analogously,
\begin{eqnarray}
\inf_x  \liminf_{n \rightarrow \infty} \int_{-n}^x \dot{g} d\lambda & = & 
\inf_x  g(x) - \limsup_{n \rightarrow \infty} g(-n) \nonumber.
\end{eqnarray}
From this and Lemma~\ref{lem:diam} we obtain
$\|\dot{g}\|_D = r < \infty$.
The uniqueness follows from the fundamental theorem of Lebesgue integral calculus and the fact that the integration constant 
is determined by the restriction  $\inf_x g(x)=0$.
 $\Box$

\section{The Dual Space $\mathscr{E}^*_{\tiny \mbox{D}}$}
\label{s:dual}
One of the central questions in functional analysis is the characterization of the dual space $\mathscr{V}^*$ of a given vector space $\mathscr{V}$.  $\mathscr{V}^*$ consists of all linear functionals $L:\mathscr{V}\rightarrow \mathbb{R}$, together with the vector space structure of pointwise addition and scalar multiplication by constants.

If the vector space $(\mathscr{V}, \|.\|)$ is equipped with a norm  $\|.\|$, 
 the question arouses about the dual norm $\|.\|^*$ in $\mathscr{V}^*$, which
is induced by 
\begin{equation}
\label{eq:operatornorm}
\|L\|^*:= \sup \{|L(\vec x)| \mid \|\vec x\|\leq 1\}.
\end{equation}
Note that $\|L\|^*$ exists if the linear functional $L$ is bounded w.r.t the norm $\|.\|$, i.e., there is a constant $M>0$ such that $|L(\vec x)| \leq M \cdot \|\vec x\|$ for all $\vec x \in \mathscr{V}$ with $\|\vec x\|$.

In this section we will determine the dual space $\mathscr{E}^*_{\tiny \mbox{D}}$ and the corresponding dual norm~(\ref{eq:operatornorm}).

First of all, consider a linear functional $L:\mathscr{E}_{\tiny \mbox{D}} \rightarrow \mathbb{R}$, 
an event sequence $\eta \in \mathscr{E}_{\tiny \mbox{D}}$ and an interval $[a,b]$.
Note that the subset of events of $\eta$ contained in $[a,b]$ defines also an event sequence.
We denote this event sequence by 
$$
\eta|_{[a,b]} := \left\{ 
\begin{array}{lcr}
\eta(t) & \ldots & t\in [a, b] \\
0       & \ldots & \mbox{else}. 
\end{array}
\right.
$$
There are only finitely many events in $[a,b]$, say at $t_{i_j} \in [a,b]$.
For convenience we write $(\eta_{t_{i_j}})_j := \eta|_{\{t_{i_j}\}}$.
So, $\eta_{t_{i_j}} \in  \mathscr{E}_{\tiny \mbox{D}}$ denotes that singleton event sequence that is zero everywhere except at $t_{i_j}$, where the event is given by $\eta({t_{i_j}})$. 
For convenience, let us write $f_L(t_{i_j}):= L(\eta_{t_{i_j}}) \in \mathbb{R}$.

The following Lemma~\ref{lem:L} is a direct consequence of the linearity of $L$.
\begin{lemma} {(\bf{Linear Functionals on $\mathscr{E}$})}
\label{lem:L}

$L$ is a linear functional on $\mathscr{E}$ if and only if there is a unique function 
$f_L:[0, \infty) \rightarrow \mathbb{R}$, such that
for all $[a,b] \subseteq [0,\infty)$ and all event sequences $\eta \in \mathscr{E}$ there holds
$$
L(\eta|_{[a,b]}) = \sum_{t_{i_j}\in [a,b]} f_L(t_{i_j}) \cdot \eta_{t_{i_j}}.
$$
\end{lemma}

Next we characterize those linear functionals which are bounded w.r.t $\|.\|_{\tiny \mbox{D}}$.

\begin{theorem} {(\bf{Bounded Linear Functionals on $\mathscr{E}_{\tiny \mbox{D}}$})}
\label{th:Lbounded}

$L \in \mathscr{E}^*$ is bounded w.r.t the discrepancy norm $\|.\|_{\tiny \mbox{D}}$
if and only if $f_L$ has bounded variation, i.e., $\|f_L\|_{\tiny \mbox{BV}} < \infty$.  
\end{theorem}

Suppose that $L$ is bounded. 
Indirectly, suppose that $\|f_L\|_{\tiny \mbox{BV}} = \infty$. 
Then there is a sequence of partitions $\mathscr{P}_k = \{t_{1}^{(k)}, \ldots, t_{n_k}^{(k)}\}$ such that
\[
\sup_{k \rightarrow \infty} \sum_{j=1}^{n_k} |f_L(t_{j+1}^{(k)}) - f_L(t_{j}^{(k)})| = \infty.
\]
This means that either (the summation is taken over all defined indexes)
\begin{equation}
\label{eq:2j}
\sup_{k \rightarrow \infty} \sum_{j} |f_L(t_{2*j}^{(k)}) - f_L(t_{2j-1}^{(k)})| = \infty
\end{equation}
or
\begin{equation}
\label{eq:2j+1}
\sup_{k \rightarrow \infty} \sum_{j} |f_L(t_{2*j+1}^{(k)}) - f_L(t_{2j}^{(k)})| = \infty.
\end{equation}

Note that  
\begin{equation}
\label{eq:modsum}
|f(t_{i+1}) - f(t_i)|=f(t_i) \eta(t_i) + f(t_{i+1}) \eta(t_{i+1}),
\end{equation}
where $(\eta(t_i), \eta(t_{i+1})) := (1,-1)$ if $f(t_i) \geq f(t_{i+1})$ and  
$(\eta(t_i), \eta(t_{i+1})) := (-1,1)$ if $f(t_i) < f(t_{i+1})$.

(\ref{eq:modsum}) together with (\ref{eq:2j}), (\ref{eq:2j+1}) means that 
there is an event sequences $\eta^{(k)}$ such that the corresponding sequence of summations
$
\psi_k := \sum_j f_L(t_j) \eta^{(k)}(t_j) = L(\eta^{(k)})
$
is unbounded, which contradicts the assumption that $L$ is bounded.
Hence, $\|f_L\|_{\tiny \mbox{BV}} < \infty$.

Now, suppose that  $\|f_L\|_{\tiny \mbox{BV}} < \infty$. 
Consider an event sequence $\eta$ with $\|\eta\|_{\tiny \mbox{D}}\leq 1$.
This means that the corresponding sequence of events $(\eta(t_i))_i$ is alternating in sign.

Consequently, we obtain
\begin{eqnarray}
\label{eq:DV}
 & & |L(\eta|_{[a,b]})| \nonumber\\
	& =&  |\sum_{t_{i_j}\in [a,b]} f_L(t_{i_j}) \cdot \eta_{t_{i_j}}| \nonumber\\
  & \leq &  |f_L(t_1) - f_L(t_2)| +  \ldots +  |f_L(t_{n-1}) - f_L(t_{n})|\nonumber\\
	 & \leq & \|f_L\|_{\tiny \mbox{BV}},
\end{eqnarray}
for any interval $[a,b]$ and any choice of partitions $(t_1, \ldots, t_n)$.
Hence, $L$ is bounded.
$\Box$

(\ref{eq:DV}) implies the following Proposition~(\ref{prop:DV}).
\begin{proposition}({\bf{Dual Discrepancy Norm}})
\label{prop:DV}

Let $f$ be of bounded variation, $\|f\|_{\tiny \mbox{BV}}< \infty$. Then
$$
\frac{1}{2} \|f\|_{\tiny \mbox{BV}}  \leq \|f\|_{\tiny \mbox{D}}^* \leq \|f\|_{\tiny \mbox{BV}},
 $$
where 
$$ 
\|f\|_{\tiny \mbox{D}}^* = \sup_{\eta \in \mathscr{E}_{\tiny \mbox{D}}, \eta\neq \vec 0}\frac{|L_f(\eta)|}{\|\eta\|_{\tiny \mbox{D}}} 
$$
and $L_f$ is the corresponding linear functional induced by $f$. 
\end{proposition}

Note that if $f$ is monotonic, we get $\|f\|_{\tiny \mbox{D}}^* = \|f\|_{\tiny \mbox{BV}}$, and if $f$ is periodically oscillating such as
$f(t) = \sin(t)$ on $[a,b]$ we obtain a low measure. This means that
\begin{equation}
\label{eq:mon}
\mu_{\tiny \mbox{mon}}(f) := \frac{\|f\|_{\tiny \mbox{D}}^*}{\|f\|_{\tiny \mbox{BV}}}
\end{equation}
measures to which extent a non-constant function $f$, i.e., $\|f\|_{\tiny \mbox{BV}}>0$,  is monotonic.
Note that  $\mu_{\tiny \mbox{mon}}(f)$ can be computed in $O(n)$ if $f$ is discrete given by $n$ values.
This can be achieved by identifying local extremal points of $f$ as points of events (``up''-event for local maximum and 
``down''-event). 
Compared to monotonicity measures~\cite{Qoyyimi14} based on rearranging the ordering in order to achieve monotonicity which is 
of $O(n\log(n))$, our measure (\ref{eq:mon}) distinguishes by its low computational complexity of $O(n)$.
A detailed study of the features of this monotonicity measure will be postponed to future research.

Is there an equivalent discrepancy measure that yields the total variation as its dual norm? 
Yes! We just need a slight modification of the discrepancy norm, the Alexiewicz norm~\cite{Alexiewicz48}:
\begin{equation}
\label{eq:AA}
\|(x_1,\ldots, x_n)\|_{\tiny \mbox{A}} := \max_k\{|\sum_{i=1}^k x_i|\}
\end{equation}
Note that $\|\vec x\|_{\tiny \mbox{D}}/2 \leq  \|\vec x\|_{\tiny \mbox{A}} \leq  \|\vec x\|_{\tiny \mbox{D}}$.

Note that $-1, 2, -2 \ldots$ is a sequence with Alexiewicz norm $1$.
Let $t_i$  mark locations of local minimum or maximum, which are alternating.
Let $t_0$ denote the first local extremum. If $t_0$ marks a local minimum, then we set $\eta(t_0):= -1$, if it is a local maximum
we set $\eta(t_0):=+1$. Then we proceed by consecutively assigning $\pm 2$ alternating in sign at the positions $t_i$ of local extrema.
By this we obtain  
\[
\sum_i f(t_i) \eta(t_i) = \sum_i |f(t_{i+1} - f(t_i))|,
\]
hence 
\begin{equation}
\label{eq:A*}
\|L_f\|_A^* = \sup_{\eta \in \mathscr{E}_{\tiny \mbox{A}}, \eta\neq \vec 0}\frac{|L_f(\eta)|}{\|\eta\|_{\tiny \mbox{A}}} = \|f\|_{\tiny \mbox{BV}},
\end{equation}
where $L_f$ denotes the linear functional induced by $f$.

\section{A Heisenberg-type Inequality between Discrepancy and Total Variation}
\label{s:inequality}

Consider $\vec x \in \{-1,1\}^n$, which is not constant, that is  $\|\vec x\|_{\tiny \mbox{BV}} >0$.
First, let us characterizes sequences $\vec x$ of minimal total variation $\|\vec x\|_{\tiny \mbox{BV}} = 2$.
$\|\vec x\|_{\tiny \mbox{BV}} = 2$ is the case if and only if there is one change in sign, hence, up to choosing the initial sign, we have
\[
\vec x\ = (\underbrace{1, \ldots, 1}_{k}, \underbrace{-1, \ldots, -1}_{n-k}).
\]
Note that $ \frac{n}{2} \leq  \max\{k, n-k\} = \|\vec x\|_{\tiny \mbox{D}}$, hence 
$n \leq \|\vec x\|_{\tiny \mbox{D}} \cdot \|\vec x\|_{\tiny \mbox{BV}}$. 
For an arbitrary number $S$ of changes in the sign we have $\|\vec x\|_{\tiny \mbox{BV}} = 2\cdot S$ and 
\[
\frac{n}{S} \leq  \max\left\{k_1, \ldots, k_S: \sum_i k_i = n, k_i \in \mathbb{N}\right\} = \|\vec x\|_{\tiny \mbox{D}}.
\]
Consequently, we obtain 
\[
n \leq \frac{n}{S} S\cdot 2 \leq \|\vec x\|_{\tiny \mbox{D}} \cdot \|\vec x\|_{\tiny \mbox{BV}}. 
\]
This result also applies to  $\vec x \in \{-1,0,1\}^n$. 
To see this, first cancel all zeros from $\vec x$ which yields $\vec {\hat x}$ and note that
\begin{equation}
\label{eq:heisenbergtype1}
\|\vec x\|_1 = \|\vec {\hat x}\|_1 \leq \|\vec {\hat x}\|_{\tiny \mbox{D}} \cdot \|\vec {x}\|_{\tiny \mbox{BV}}
\leq \|\vec x\|_{\tiny \mbox{D}} \cdot \|\vec x\|_{\tiny \mbox{BV}}.
\end{equation}

Finally, (\ref{eq:heisenbergtype1}) implies the Heisenberg-type inequality between discrepancy and total variation
\begin{equation}
\label{eq:heisenbergtype}
\|\vec x\|_1 \leq \|\vec x\|_{\tiny \mbox{D}} \|\vec x\|_{\tiny \mbox{BV}}
\end{equation}
for event sequences $\vec x \in \{-1,0,1\}^n$, $n \in \mathbb{N}$ and $\|\vec x\|_{\tiny \mbox{BV}}>0$.

This inequality is sharp as  $\vec x = (1, \ldots, 1,0,\dots,0)$ or 
$\vec x = (1, \ldots, 1,-1,\dots,-1)$ induce equality.
By interpreting $\|\vec x\|_{\tiny \mbox{BV}}$ in (\ref{eq:heisenbergtype}) as measure of oscillation and 
 $\|\vec x\|_{\tiny \mbox{D}}$ as measure of the dominance of monolithic blocks, we recognize the similarity to 
Heisenberg's inequality relation.


\section{Conclusion}
In this article we investigated the space of event sequences as normed space equipped with 
Weyl's discrepancy norm, which distinguishes by its property to turn threshold-based sampling into a quasi-isometry mapping.
As result we found various characterizations and interpretations of this norm as for example a 
Jordan-like decomposition law.
We also investigated its relationship to total variation and found a Heisenberg-type inequality.
The ratio of the dual discrepancy norm and total variation turns out to be a measure of monotonicity, which will be investigated in more detail in the future.
 

\section*{Acknowledgment}
The author would like to thank the Austrian COMET Program and in particular Florian Sobieczky for careful reviewing and 
fruitful discussions.



%

\end{document}